\documentclass{smfart}
\usepackage{AW-38}
\begin{document}

\title[La guerre des recensions]{La guerre des recensions\\
autour d'une note d'André Weil en 1940}
\author{Mich\`ele Audin}
\address{Institut de Recherche Math\'ematique Avanc\'ee\\
Universit\'e de Strasbourg et CNRS\\
7 rue Ren\'e Descartes\\
67084 Strasbourg\index{Strasbourg} cedex\\
France}
\email{michele.audin@math.unistra.fr}
\urladdr{http://www-irma.u-strasbg.fr/~maudin}
\thanks{Version pr\'eliminaire du \today}


\newcommand\JFM{\emph{JFM}\xspace}
\newcommand\MR{\emph{MR}\xspace}
\newcommand\Zbl{\emph{Zbl}\xspace}

\begin{abstract}
Nous étudions avec précision les recensions d'une note d'André Weil de 1940 dans les trois journaux de référence \emph{Jahrbuch über die Fortschritte der Mathematik},  \emph{Zentralblatt für Mathematik und ihre Grenzgebiete} et \emph{Mathematical Reviews}, avec le contexte de leurs publications.
\end{abstract}

\begin{altabstract}
We investigate the reviews of a \emph{Comptes rendus} note by André Weil in 1940 in the three journals \emph{Jahrbuch über die Fortschritte der Mathematik},  \emph{Zentralblatt für Mathematik und ihre Grenzgebiete} et \emph{Mathematical Reviews}, together with the context of the publication of these reviews.
\end{altabstract}

\maketitle

Cet article tourne autour d'une note publiée aux \emph{Comptes rendus de l'Académie des sciences} par André Weil en 1940. Il y est bien entendu question du contenu de cette note --- l'annonce, entre autres résultats, d'une démonstration de l'\og hypothèse de Riemann pour les corps finis\fg. Mais le sujet de l'article est la description du contexte dans lequel cette note a été lue et recensée dans les trois journaux de référence qui existaient à l'époque, le \og vieux\fg \emph{Jahrbuch über die Fortschritte der Mathematik}, qui vivait ses derniers jours, et les plus modernes, allemand \emph{Zentralblatt für Mathematik und ihre Grenzgebiete} et américain \emph{Mathematical Reviews} (dont les titres sont désormais abrégés en  \JFM, \Zbl, \MR). 

\medskip
On le sait, les recensions ne sont pas isolées de la situation politique et sociale. Le rôle qu'a joué le contexte géopolitique dans la création de \emph{MR} aux \'Etats-Unis en 1940 après celle de \Zbl en Allemagne en 1931 a été particulièrement étudié (voir~\cite{SiegmundSBerichter,SiegmundS94}). Nous ne reviendrons pas sur cette histoire, mais le présent article pourra être considéré comme une sorte d'\og étude de cas\fg (un cas un peu extrême) illustrant cette question. C'est la rivalité entre les deux mathématiciens Helmut Hasse et André Weil, l'un allemand, officier de marine occupé à travailler pour l'armée à Berlin, l'autre français et insoumis (n'ayant pas répondu à l'ordre de mobilisation), puis juif réfugié aux \'Etats-Unis, et peut-être aussi de deux écoles, l'une plus algébrique, l'autre plus géométrique, dans ce qu'il faut bien appeler une \og guerre des recensions\fg, qui est envisagée ici.

\medskip
Les deux auteurs avaient compris l'importance de ces journaux. Notons que Hasse était impliqué dans la recherche de recenseurs pour le \Zbl. Weil, nous allons le voir assez en détail, les a utilisés, comme lecteur, pour leur contenu (en 1940) puis comme auteur, dès son arrivée aux \'Etats-Unis, pour avoir connaissance rapidement des nouvelles publications, voire donner son avis.

\medskip
De cette note fameuse d'André Weil (et d'au moins une, d'ailleurs fameuse elle aussi, de ses recensions) il a déjà été beaucoup question. Outre le commentaire de l'auteur lui-même dans ses \OE uvres complètes~\cite{weilOC1}, mentionnons ici seulement l'article~\cite{schappa06}. L'ambition ici est simplement d'éclairer un peu le paysage de ces recensions, en mettant la chronologie en évidence d'abord, en utilisant cette chronologie ensuite, et en ajoutant des informations issues des correspondances des protagonistes, non seulement sur les recensions et leurs contextes, mais aussi sur la façon dont eux-mêmes utilisaient les journaux d'analyse.

\medskip
Les sources utilisées pour cet article sont
\begin{itemize}
\item les textes publiés d'André Weil~\cite{weil40,weilOC1,weilsouvenirs} et de Helmut Hasse~\cite{Hasseoslo,Hasse42}, 
\item les recensions publiées des articles~\cite{weil40,Weil41,Hasse42,Hasse42a,Weil48a,weilOC1}, que j'ai lues dans les volumes papier (à la bibliothèque de l'\textsc{irma} de Strasbourg),
\item la correspondance publiée d'André Weil avec Henri Cartan~\cite{CW-doc} et 
\item la correspondance de Helmut Hasse\footnote{Précisons que les sources utilisées sont les copies carbone, conservés par Hasse, des lettres dactylographiées expédiées à ses correspondants.} avec Gaston Julia, André Weil, Bartel van der Waerden, le Rektor de son université (Nachla\ss\ Hasse, Université de Göttingen) et avec Wilhelm Süss (Nachla\ss\ Süss, Université de Freiburg)
\end{itemize}
ainsi que diverses sources publiées, qui sont mentionnées, au fur et à mesure (avec renvois à la bibliographie).

\begin{remarque*}
Ayant mentionné les versions \og papier\fg, je me permets de faire quelques remarques sur les versions numérisées des journaux de recension:
\begin{itemize}
\item La numérisation du \JFM ne semble donner aucune indication sur la date de la recension (il n'y en a pas non plus dans le volume papier, mais il y a au moins une date de publication de l'ensemble).
\item Celle de \Zbl donne accès à la recension dans le volume papier scanné (tout entier, on peut regarder la page précédente, la suivante, feuilleter), ce qui est extrêmement intelligent et pratique.
\item C'est après beaucoup de temps passé sur la version en ligne de \MR que j'ai fini par comprendre que l'information sur la date de la publication était bien donnée. Pour faire gagner du temps aux éventuels lecteurs de cet article, par exemple, dans

\begin{quote}
{\bfseries MR0004766 (3,58e)}

Ehresmann, Charles; Feldbau, Jacques

{\bfseries Sur les propriétés d'homotopie des espaces fibrés. (French)}

C. R. Acad. Sci. Paris 212, (1941). 945--948. (Reviewer: A. Weil)

56.0X 
\end{quote}
l'information voulue est dans le 3,58e, qui signifie que la recension est parue dans le volume 3, à la page 58, et qu'elle est la cinquième (e) de cette page. Sachant que 1 = 1940, on en déduit que 3 = 1942. Je n'avais pas trouvé comment demander au moteur de recherche de trouver les recensions qui sont sur la même page, mais, comme me l'a indiqué un lecteur d'une première version de ce texte, entrer le numéro du volume et celui de la page (ici \og 3,58\fg) dans le champ \og numéro du rapport\fg fournit bien la liste voulue. Hélas, l'affiliation du recenseur (dans cet exemple, A. Weil (Haverford, Pa.)) est perdue. 
\end{itemize}
\end{remarque*}
$$\star$$

\medskip
On ne trouvera ici de biographie ni de l'un, ni de l'autre, des deux principaux protagonistes. En particulier, rien ne sera dit, ni de la carrière de Helmut Hasse, ni de son rôle de directeur de l'institut de mathématiques de  Göttingen à partir de 1935, ni même de sa participation à la crise du \Zbl en 1938, sujets pour lequels on pourra consulter~\cite{SchappaGottingen,SiegmundSBerichter}. La question de ses choix politiques ne sera pas discutée non plus. La figure du mathématicien aux prises avec la guerre dans laquelle il est enrôlé sera prédominante (nous le verrons rappelé dans la \emph{Wehrmacht} en avril 1940, en mission scientifique dans Paris occupé à l'automne 1940, employé, à Berlin, à d'autres tâches mathématiques que la théorie des nombres ensuite). 

De même, les raisons qui ont poussé André Weil à refuser de rentrer en France au moment de la mobilisation de septembre 1939 et les circonstances qui ont fait qu'il s'est quand même retrouvé dans une prison française, puis, par la vertu de la législation de Vichy, aux \'Etats-Unis (voir~\cite{weilsouvenirs,weilOC1,Pekonen,CW-doc,AW-38}) ne seront pas détaillées ici.

Les aspects \og italiens\fg (mathématiques et/ou politiques) de cette histoire ne seront pas envisagés non plus. Voir encore~\cite{schappa06} pour le rôle de Severi (dont le nom apparaîtra bien sûr ici ou là dans ce texte, \emph{via} ses travaux mathématiques).

\medskip
Le contexte mathématique (et même humain) tel qu'il était peu d'années auparavant (et  avant la guerre), en 1936, est résumé par Norbert Schappacher dans~\cite{schappa06}. Résumons. En 1932, Hasse a démontré un théorème analogue à l'hypothèse de Riemann, pour la fonction $\zeta$ attachée au corps des fonctions sur une courbe elliptique sur un corps fini. L'histoire de ce théorème et des démonstrations de Hasse est le sujet de trois articles très détaillés de Roquette~\cite{Roquette02,Roquette04,Roquette06}, auxquels nous renvoyons. La question restait ouverte de démontrer ce résultat pour les courbes de genre supérieur\footnote{Les fonctions $\zeta$ généralisées remontent à la thèse d'Emil Artin en 1924. Nous suivons ici la description de Hasse~\cite{Hasseoslo}. Un corps fini $k=\FF_q$ et une courbe algébrique lisse sur $k$ sont fixés, $N_n$ désigne le nombre de points fermés de degré $n$ de la courbe, la fonction $\zeta$ est définie comme le produit eulérien
$$\zeta(s)=\prod_{n=1}^{\infty}\lefpar\frac{1}{1-\dfrac{1}{q^{ns}}}\rigpar^{N_n}.$$
Le cas de genre $0$, où la courbe est la droite projective, donne $N_n=q^n+1$ et
$$\zeta_0(s)=\frac{1}{1-\dfrac{1}{q^s}}\frac{1}{1-\dfrac{1}{q^{s-1}}}$$
et, grâce au théorème de Riemann-Roch, 
$$\zeta(s)=\zeta_0(s)P\lefpar\dfrac{1}{q^s}\rigpar$$
pour un certain polynôme $P$ de degré $2g$ ($g$ est le genre de la courbe). Les zéros de $\zeta$ sont les $s$ tels que $q^{-s}$ soit une des racines de $P$. L'hypothèse de Riemann généralisée postule que la partie réelle de ces zéros est $1/2$, c'est-à-dire que les racines de $P$ sont de module $q^{-1/2}$.}. Une courbe elliptique possède un anneau d'endomorphismes, que les courbes de genre plus grand n'ont pas et qui jouait un rôle important dans la démonstration de Hasse. Il fallait donc un analogue des endomorphismes en genre plus grand. Pour ceci, Max Deuring (1907--1984) avait proposé d'utiliser la théorie des correspondances sur ces courbes. Weil était au courant de tout ceci, il avait parlé avec ces mathématiciens allemands, il correspondait avec Hasse et, d'ailleurs, ce dernier lui avait écrit, le 12 juillet 1936, une longue lettre (huit pages!) dans laquelle il expliquait ce programme\footnote{\label{notedoubles}Rappelons que c'est la copie conservée par Hasse que nous lisons. Il se peut que \og le côté Weil\fg de cette correspondance  ait été égaré au cours des déplacements de Weil avant et après 1940. D'autre part, Weil a beaucoup jeté. En tout cas, nous ne savons pas où se trouvent les lettres qu'il a reçues de Hasse. Signalons que, même dans le Nachla\ss\ Hasse, il manque des lettres. Sur l'utilisation de cette lettre par Weil, voir plus bas page~\pageref{weiletlalettre}.}, que Hasse lui-même mentionnait en conclusion de son exposé~\cite{Hasseoslo} au congrès international d'Oslo (auquel Weil ne participa pas) quelque jours plus tard. Hasse organisa un petit colloque sur la géométrie algébrique en janvier 1937 à Göttingen, avec notamment Deuring, van der Waerden et Geppert (mais sans Weil\footnote{Nous ne savons pas si Weil avait vraiment envisagé de s'y rendre. Il écrivit à Hasse le 21 décembre 1936:

\begin{quote}
{\smaller Ich hätte gern vesucht, nach Göttingen zum Kolloquium über algebraische Funktionen zu fahren (wo eben, wie es scheint, Deuring diese Resultaten vortragen soll) wenn ich mich nicht eben am 6. Januar nach Amerika einschiffen müsste.

}
\end{quote}
mais il dit aussi dans ses souvenirs qu'il avait renoncé, dès 1934, à se rendre dans l'Allemagne nazie.\label{notegottingen}}), il exposa à nouveau ses travaux pendant la réunion de la \textsc{dmv} à Baden-Baden en septembre 1938, puis au cours de trois exposés à Paris (où Gaston Julia l'avait invité) en mai 1939, sous le titre général (en français) \og Nouvelles recherches sur l'arithmétique des corps de fonctions algébriques\fg. Les titres des trois exposés étaient\footnote{On trouve ces titres d'exposés dans le rapport que Hasse a expédié à son Rektor et qui est conservé dans son Nachla\ss.} (en français\footnote{Toutes les citations de ce texte sont données dans leur langue originelle.})
\begin{enumerate}
\item Généralités: Le groupe des classes de diviseurs et l'anneau des multiplicateurs.
\item Points rationnels et entiers sur les courbes algébriques à coefficients entiers.
\item Points rationnels sur les courbes algébriques à coefficients mod $p$.
\end{enumerate}
Mais André Weil était déjà parti pour la mission qui devait le conduire jusqu'au nord de la Finlande, puis jusqu'à une prison française en 1940. Ce qui nous amène, nous aussi, en avril 1940.

\subsection*{Avril 1940 (la note de Weil)}
En avril 1940, André Weil était incarcéré à la prison de Bonne-Nouvelle à Rouen et il y attendait son procès pour insoumission.  Il était approvisionné en références bibliographiques par les visites de sa famille (parents, s\oe ur, épouse) et par la correspondance avec ses amis mathématiciens et notamment avec Henri Cartan. Puisqu'il est question ici de recensions, signalons que, comme le montre la correspondance, André Weil disposait, dans sa prison, du tout premier fascicule de \MR, qui venait de paraître (et était bien arrivé en France: Cartan écrivit le 12 février qu'il venait de le recevoir). D'autre part, Henri Cartan lisait le \Zbl pour lui, à sa demande. Il n'est pas hors-sujet de citer quelques passages de leurs échanges. Le 17 mars, par exemple, Weil écrivit à Cartan
\begin{quote}
Peux-tu puiser pour moi, dans le Journal de Crelle, le renseignement suivant: quel est, dans le groupe des classes de diviseurs sur une courbe algébrique (pour un corps de constantes \emph{fini} à $q=p^n$ éléments) le nombre de classes dont la $m$-ième puissance est l'unité? Je ne suis pas sûr que le problème ait été résolu pour un genre quelconque, mais il l'est en tout cas pour le genre~$1$, soit par Hasse lui-même, soit par un de ses élèves. Il n'y a qu'à voir les titres des mémoires parus dans Crelle depuis 5 ou 6 ans.
\end{quote}
Et Cartan lui répondit, de Clermont-Ferrand, le 2 avril:
\begin{quote}
Heureusement nous avons le Zentralblatt jusqu'au tome 15 (1937) inclus\footnote{Cartan parlait ici de la bibliothèque de l'université de Strasbourg, qui était repliée à Clermont-Ferrand depuis septembre 1939, et signalait qu'une partie des livres étaient restés à Strasbourg chez un relieur. Une histoire passionnante mais un peu éloignée de notre sujet.} [...] Voici ce que j'ai trouvé dans le Zbl:
\begin{enumerate}
\item \emph{mémoire de Hasse} (Zur Theorie des abstrakten ellipt. Funktionenkörper; Crelle, 175, p. 55--62; Zbl, 14, p. 149) Résultat: pour un corps $k$ de caractéristique $p$, algébriqu\up{t} fermé, et un corps $K$ elliptique de fonctions, le nombre des classes de diviseurs dont la $n$\up{e} puissance est l'unité est $$\begin{array}{ll} n^2&\text{ si }n\not\equiv 0(p)\\ n&\text{ si }n=p^\nu\text{ et }A\neq 0~;\\ 1&\text{ si }n=p^\nu\text{ et }A=0~;
\end{array}$$
($A$ est défini par le développement d'un élément $v\in K$, $v=\dfrac{1}{\pi^p}-\dfrac{A}{\pi}+\cdots$)
\item \emph{mémoire de Deuring} (Automorphismen und Divisorenklassen der Ordnung $\ell$ in algebraischen Funkt.-körpern; \emph{Math. Ann.} 113, 1936, p. 208-215). Le compte-rendu (Zbl 14, p. 293) est fait par Hasse\footnote{Ce qu'a écrit Cartan sur le mémoire de Deuring est plus ou moins une traduction de la recension de Hasse pour le \Zbl.} Il semble qu'il n'y ait aucun résultat définitif. Hasse y parle de l'\emph{hypothèse} suivante, concernant le nombre $T_p$ des solutions de $X^p=1$ pour $p$ premier: 
$$\begin{cases} T_p=p^{2g}&(g=\text{ genre}) \text{ si la caractéristique du corps }k\text{ est }\neq p;\\ T_p=p^\rho,&\text{ avec }0\leq\rho\leq g\text{ si la caract. est }p.
\end{cases}$$
Cette hypothèse est en accord, dit Hasse, avec ce que démontre Deuring: soit $K$ (corps de fonctions sur $k$ algébriq\up{t} fermé) séparable, cyclique, de degré 1\up{er} $p$ sur~$K_0$; supposons que $t$  diviseurs 1\up{ers} soient ramifiés dans $K$ ($t>0$); et que, dans $K_0$, chaque classe de diviseurs de degré $0$ soit la $p$\up{e} puissance de $p^{\rho_0}$ classes exactement. Alors: si $k$ n'est pas de caract. $p$, chaque classe de degré $0$ de $K$ est la $p$\up{e} puissance de $$p^{p\rho_0+(p-1)(t-2)}\text{ classes };$$
si $k$ est de caract. $p$, chaque classe de degré $0$ de $K$ est la $p$\up{e} puissance de $$p^{p\rho_0+(p-1)(t-1)}\text{ classes}.$$
Pour la démonstration, on se sert du \og Normensatz\fg de Tsen. Il paraît d'autre part que, dans une addition à la correction, Deuring affirme sans démonstration que la résolution de $X^n=C$ est facile à traiter dans le cas général par la méthode de Hasse. Hasse renvoie d'autre part à un travail de \emph{Hasse et Witt} paru dans Mh. Math. Phys., et initulé \emph{Zyklische unverzweigte Erweiterungs\-körper vom Primzahlgrade $p$ über einem alg. F.-körper des Char. $p$}. 
\end{enumerate}

Voilà. C'est tout ce que j'ai trouvé. J'espère que ce n'est pas en dehors du sujet. J'attends maintenant tes commentaires; dois-je tenter autre chose? Ne pourrais-tu te faire apporter de Paris les numéros du Zbl à partir du tome 16 (inclus)?
\end{quote}
\`A quoi Weil répondit, le 5 avril:
\begin{quote}
Merci de ta lettre et du mal que tu t'es donné. Le résultat de Hasse est bien celui que je voulais avoir, pour le genre~$1$. Pour le genre $>1$, je suis à peu près certain que le mémoire de Deuring annoncé d'après ce que tu m'écris, en 1936, a dû paraître\footnote{Il est paru en 1937 (avec des corrections pour lesquelles Deuring a remercié Hasse et Schmid). Il s'agit de~\cite{Deuring37}. Mais le résultat voulu par Weil n'y figurait pas. Dans la lettre à Hasse datée de décembre 1936 (et déjà citée dans la note~\ref{notegottingen}), que Weil écrivit juste avant de partir pour un séjour à Princeton en 1937, il avait demandé si les résultats de Deuring allaient paraître rapidement.}. Mais comment faire pour le savoir? Je suppose que tu iras à Paris un de ces jours et pourras me renseigner. Me procurer le Zentralblatt ici est presque impraticable. En attendant, je vais admettre le résultat ou bien le démontrer (je n'ai besoin que du cas \og trivial\fg de $n$ premier à la caractéristique du corps). Mon travail fait des progrès surprenants. A peu près toute la théorie \emph{transcendante} des fonctions algébriques se transpose aux corps de constantes finis, avec matrices des périodes, relations bilinéaires, et le théorème de Hurwitz; la théorie algébrique de Severi marche trivialement, ce qui n'a pas empêché Deuring, qui, comme tous ces gens-là ne connaît pas ses classiques, de la retrouver en son langage, non sans grands cassements de tête, idéaux, diviseurs, corps de classes de restes, etc., il y a 3 ou 4 ans).
\end{quote}
Le fond mathématique du débat, algèbre \emph{vs} géométrie algébrique, apparaît clairement ici. Le 8 avril, Weil avait rédigé la fameuse note et il écrivit à Cartan:
\begin{quote}
L'absence de renseignements précis sur les travaux de Hasse et Deuring (nombre de $n$\up{ièmes} de périodes et même de racines de $x^n=a$) me gêne fort\footnote{Nous l'avons vu, Cartan, à la suite de Hasse, avait écrit qu'il n'y avait aucun résultat définitif sur les points de $n$-division dans l'article~\cite{Deuring36}. Mais visiblement, André Weil pensait que Hasse et Deuring avaient pu démontrer ultérieurement que l'ordre des points de $n$-division (dans la jacobienne d'une courbe de genre $g$) est $n^{2g}$. Comme me l'a expliqué Jean-Pierre Serre lorsque je finissais d'écrire les notes de~\cite{CW-doc}, Weil se rendit compte plus tard qu'il n'en était rien et qu'il s'agit d'un théorème difficile et important. La lettre qu'il écrivit à Artin (voir~\cite{weilOC1}) en 1942 montre que, deux ans plus tard, il savait le démontrer.}; il me faudrait, non seulement le résultat (pour $n$ premier à la caractéristique, du moins) mais des indications sur la méthode. Penses-tu aller à Paris prochainement ou y aurait-il quelque normalien assez dégourdi pour faire le nécessaire? provisoirement, rédigeant une note aux Comptes Rendus que je viens d'expédier à ton père\footnote{Parce qu'Élie Cartan (1869--1951) était membre de l'Académie des sciences.}, j'ai pris le parti d'admettre tout bonnement le résultat, et sans le dire; mais pour la suite de mon travail je crois qu'il me serait utile d'avoir des précisions sur la démonstration, à part le fait qu'il serait tout de même bon d'en être sûr, puisque tout ce que je fais repose là-dessus. Chose plus sérieuse, j'ai expédié la note sans attendre d'avoir démontré le lemme fondamental; mais j'y vois assez clair à présent sur ces questions pour en prendre le risque. Jamais je n'ai rien écrit, et je n'ai presque jamais rien vu, qui atteigne un aussi haut degré de concentration que cette note. Hasse n'a plus qu'à se pendre, car j'y résous (sous réserve de mon lemme) \emph{tous} les principaux problèmes de la théorie: 1) hypothèse de Riemann pour les fonctions $\zeta$ de ces corps (démontrée par Hasse pour le genre~$1$); 2) les séries $L$ d'Artin relatives aux caractères des extensions algébriques de ces corps sont des polynomes, dont je détermine le degré (connu, en vertu de la loi de réciprocité d'Artin, pour les extensions \emph{abéliennes} seulement). Le second point est le plus intéressant pour moi, car il ouvre la voie pour l'étude des extensions non abéliennes de ces corps, et, en vertu des analogies très étroites avec l'arithmétique, des corps de nombres algébriques. Mais le premier point est plus sensationnel, à cause de l'acharnement qu'ont mis Hasse et tout son monde à le démontrer. Ils en sont passés assez près, sans le voir; ce qui leur a manqué, c'est l'analogue pour ces corps de la théorie transcendante des correspondances.
\end{quote}

On le voit, Weil savait très bien que son lemme n'était pas démontré. Mais il avait confiance dans son approche. Comme il l'écrivit dans ses \OE uvres complètes (à un moment, bien sûr, où il savait que cette approche avait porté ses fruits):
\begin{quote}
\noindent En d'autres circonstances, une publication m'aurait semblé bien prématurée. Mais en avril 1940, pouvait-on se croire assuré du lendemain?~\cite[p.~546]{weilOC1}\footnote{Le paragraphe complet dans lequel cette phrase a été prélevée est cité à la fin du présent article.}
\end{quote}

La note a été transmise par \'Elie Cartan à l'Académie des sciences et elle est parue dans le \emph{Compte rendu} daté du 22 avril 1940. 

\medskip
Pendant ce temps, Hasse devait rejoindre son poste dans la Wehrmacht:
\begin{quote}
Als ich im April 1940 zum Wehrmachtsdienst einberufen wurde [...]~\cite{Hasse42}\footnote{Le paragraphe complet dans lequel cette phrase a été prélevée est cité un peu plus bas (page~\pageref{Wehrmachtsdienst}).}
\end{quote}
Nous le verrons plus tard reprocher à Weil d'avoir \og profité\fg de ce moment et de la guerre pour publier cette note. 

\subsection*{Octobre-novembre 1940 (réactions de Hasse)}
En octobre et novembre 1940, Hasse est venu deux fois à Paris occupé (dans son uniforme d'officier de marine) et nous savons qu'un des buts de ces voyages était de recruter des recenseurs français pour le \Zbl\footnote{Harald Geppert (1902--1945), qui dirigeait désormais à la fois le \Zbl et le \JFM, est lui aussi venu à Paris, dans le même but, en décembre. Il avait pris la direction des journaux allemands en 1938 après le départ d'une bonne partie de la rédaction de \Zbl et de son rédacteur en chef Otto Neugebauer. Ce géomètre écrivit lui-même des centaines de recensions. C'était un nazi convaincu et il se suicida en 1945.}. Nous savons qu'il a vu, entre autres mathématiciens, Gaston Julia en octobre et Jean Dieudonné en novembre. C'est probablement lors de la première de ces deux visites qu'il a eu connaissance de la note de Weil. Si, et quand, le \emph{Compte rendu} du 22 avril est arrivé aux abonnés allemands reste à déterminer, mais nous savons (voir la lettre de septembre 1941 citée ci-dessous) que c'est Julia qui a \og procuré\fg la note à Hasse. Nous savons aussi que Hasse a passé une soirée avec Dieudonné en novembre (il l'a écrit dans une lettre adressée à Wilhelm Süss, dans laquelle il racontait son voyage, le 19 décembre). C'est sans doute à ce moment-là qu'il lui a dit ce que Weil rapportait avec satisfaction, presque quarante ans plus tard:
\begin{quote}
\noindent Si j'en crois Dieudonné, qui vit Hasse à Paris quelques mois plus tard, celui-ci s'en montra fort scandalisé.~\cite[p.~546]{weilOC1}
\end{quote}

\medskip
Pendant ce temps, André Weil, après procès, incorporation, guerre, armistice et démobilisation, était arrivé à Clermont-Ferrand, où Dieudonné le rejoignit d'ailleurs assez vite (il était présent lors du congrès Bourbaki qui se tint du 7 au 10 décembre) et où ils habitèrent le même hôtel~\cite[p.~178]{weilsouvenirs}, ce qui leur permit sans doute d'échanger quelques informations et commentaires.

\subsection*{Février 1941 (recension de la note de Weil par le \Zbl)}
C'est dans le fascicule de février 1941 de \Zbl que parut la première recension de la note. Le rapporteur était van der Waerden. Il plaça la question dans son contexte mathématique (celui des correspondances et des travaux de Hasse et Deuring), dit que la note était une annonce (une \og esquisse\fg), il énonça le lemme en question et écrivit que l'hypothèse de Riemann en découlait. Voici le texte complet.

\begin{quote}
Skizze einer Lösung der Hauptprobleme des Theorie des algebraischen Funktionenkörper mit endlichem Konstantenkörper. Wie Hasse und Deuring erkannt habe, gibt die Theorie der algebraischen Korrespondenzen den Schlüssel zu diesen Problemen, aber die Severische algebraische Theorie der Korrespondenzen genügt nicht, sondern man mu\ss\ die Hurwitzsche transzendente Theorie auf diese Funktionenkörper übertragen. Zu dem Zweck wird die Gruppe der Divisorenklassen nullter Ordnung isomorph abgebildet auf eine Gruppe von Vektoren, d.h. von einspaltigen Matrices mit $2g$ Reihen, deren Elemente einem Ring $\gotR$ angehören, modulo $1$. $\gotR$ entsteht aus dem Körper der rationalen Zahlen durch Komplettierung in bezug auf alle $l$-adischen Bewertungen, wobei $l$ alle Primzahlen mit Ausnahme der Körpercharakteristik $p$ durchläuft, und $g$ ist das Geschlecht des Funktionnenkörpers. Nunmehr kann man genau nach Hurwitz jeder algebraischen $(m_1,m_2)$-Korrespondenz eine $2g$-reihige quadratische Matrix $L$ zuordnen, deren Elemente $\equiv 0$ (mod $1$) in $\gotR$ sind. Die inverse Korrespondenz definiert ebenso eine Matrix $L'$, und unter gewissen Bedingungen ist die Spur von $LL'$ gleich $2m_2$. Mit Hilfe dieses Lemmas folgt die \og Riemannsche Vermutung\fg für Funktionenkörper. Schlie\ss lich wird die Wirkung der Galoisschen Gruppe eines solchen Körpers auf die den Divisorenklassen zugeordneten Vektoren untersucht.
\end{quote}

Si l'on doit compter les antagonismes algèbre \emph{vs} géométrie algébrique parmi les éléments de cette histoire, alors il faut signaler que Bartel Leendert van der Waerden (1903--1996) se plaçait, parmi les \og Allemands\fg, plus du côté de la géométrie algébrique que l'entourage de Hasse. Sur cette question, il n'en sera pas dit beaucoup plus, voir~\cite{schappa07,RSSvanderWaerden}. Voir toutefois aussi page~\pageref{demvdW}.


\subsection*{Avril 1941 (recension de la note de Weil par \MR)}
C'est juste après l'arrivée de Weil sur le territoire américain que \MR publia sa recension de la note. L'auteur était Schilling. Voici son texte:
\begin{quote}
The author sketches in this note a proof of Riemann's hypothesis for function fields $K=k(x,y)$ of one variable and genus $g$ over a finite field $k$ of $q$ elements. For the proof, it is necessary to generalize the theory of fixed points of algebraic correspondences $C$ on $\ov{K}=\ov{k}K$, where $\ov{k}$ denotes the algebraic closure of $k$. The multiplicative group of all $\omega\neq 0$ in $\ov{k}$ is isomorphic to the additive group of all elements mod 1 in the universal ring $\gotR$ over the field of rationals with respect to all primes different from the characteristic of $k$. By means of such a fixed isomorphism and the representation of divisors in $\ov{K}$ by residue classes of functions in $\ov{K}$, the group $G$ of all classes in $\ov{K}$ whose orders are prime to $q$ can be mapped upon the additive group of all $2g$-dimensional vectors with components in $\gotR$ mod 1. The automorphism of $\ov{K}\vert K$ which extends the Frobenius automorphism $\omega\to\omega^q$ induces a linear transformation $I$ on the vector group. The transformation $I$ can be interpreted as a $2g$-dimensional matrix with coefficients in $\gotR$. A $(m_1,m_2)$ correspondence $C$ of $K$ induces a homomorphism on $G$ which can be described by a matrix $L$ operating on the representation of $G$ by vectors. The number of fixed points of $C$ is ``in general'' equal to $m_1+m_2-\tr(L)$. The key to the author's investigation is the fact that the trace $\tr(LL')$ is a positive rational integer, where $L'$ denotes the matrix of the inverse correspondence of $C$. This result is a generalization of the fundamental relation for the periods of abelian integrals of first kind on an algebraic Riemann surface. A careful analysis of the ring of correspondences on $\ov{K}$ shows that the latter is a finite hypercomplex order over the ring of rational integers. The actual proof of Riemann's hypothesis consists of the investigation of the correspondence $\Sigma= (x\to x^q, y\to y^q)$. The matrix representing this correspondence is equal to $I$ and its characteristic polynomial $\module{E-uI}$ coincides with the polynomial factor of the zeta function of $K$. The fundamental inequality $\tr(LL')>0$ applied to the correspondences $a_0+a_1\Sigma+\cdots+a_{2g-1}\Sigma^{2g-1}$, $a_i$ rational integers, proves that all the characteristic roots have $\sqrt{q}$ as absolute value. Finally, the author mentions that his methods prove Artin's series for a finite extension of $K$ to be polynomials. The reviewer hopes that a detailed account of these important results will soon be accessible.
\end{quote}

Otto Franz Georg Schilling (1911--1973), avait été un étudiant d'Emmy Noether et Helmut Hasse, avec qui il avait obtenu sa thèse à Marburg en 1934. Il vivait aux \'Etats-Unis depuis 1935. Un deuxième algébriste allemand (si le Néerlandais van der Waerden peut être qualifié comme tel) lisait la note, constatait que c'était une annonce (le verbe anglais \emph{sketches} de Schilling était strictement équivalent au substantif allemand \emph{Skizze} de van der Waerden) de résultats importants et le disait. 

Avant d'arriver à la troisième recension, il importe de suivre la chronologie.

\subsection*{Juin-juillet 1941 (la note \og américaine\fg de Weil)}
Arrivé aux \'Etats-Unis grâce à la mission Rapkine (voir~\cite{AW-38}) et installé provisoirement (dans le bureau de Claude Chevalley) à Princeton, André Weil s'était remis au travail. La note~\cite{Weil41} qu'il présenta à la National Academy le 1\up{er} juin commence ainsi:
\begin{quote}
A year ago I sketched the outline of a new theory of algebraic functions of one variable over a finite field of constants, which may suitably be described as transcendental, in view of its close analogy with that portion of the classical theory of algebraic curves which depends upon the use of Abelian integrals of the first kind and of Jacobi's inversion theorem; and I indicated how this led to the solution of two outstanding problems, viz., the proof of the Riemann hypothesis for such fields, and the proof that Artin's non-abelian $L$-functions on such fields are polynomials. I have now found that my proof of the two last-mentioned results is independent of this ``transcendental'' theory, and depends only upon the 
algebraic theory of correspondences on algebraic curves, as due to Severi.
\end{quote}
Le fascicule des \emph{Proceedings} de la {National Academy} dans lequel cette note est parue est daté du 15 juillet 1941. Dieudonné connaissait la nouvelle et en avait fait part à Cartan\footnote{Si la correspondance entre Cartan et Weil a complètement cessé entre le départ de Weil aux \'Etats-Unis et la Libération de Paris, il semble que Weil ait pu communiquer avec Dieudonné, peut-être par un intermédiaire, et peut-être parce que Dieudonné était en zone non occupée alors que Cartan était à Paris. Les désignations \og André\fg et \og H.H.\fg évoquent une possible censure: les amis bourbakistes ne s'appelaient pas, en temps normal, par leurs prénoms.} dans une lettre du 7 juillet 1941:

\begin{quote}
J'ai eu récemment des nouvelles d'André; il a enfin trouvé une démonstration complète de sa Note de l'an dernier, il arrive paraît-il à se passer du fameux lemme (qui, lui, reste toujours en suspens); je suis ravi de cette nouvelle, particulièrement parce que ça va faire crever H.H. de jalousie! [...]
\end{quote}

\subsection*{Septembre 1941 (réactions de Hasse)}
La publication ne dut pas être beaucoup plus tardive puisque, en septembre, Hasse avait déjà reçu un exemplaire de la nouvelle note, cette fois envoyée par Weil lui-même, avec une lettre d'accompagnement. Nous le savons parce qu'il l'écrivit à son ami Gaston Julia le 14 septembre\footnote{Lettre à Gaston Julia, écrite les 7 et 14 septembre 1941, Nachla\ss\ Hasse, Université de Göttingen. Comme dans le cas de la correspondance avec Weil (voir la note~\ref{notedoubles}), il s'agit du double conservé par Hasse. Nous ne savons pas si \og le côté Julia\fg a été conservé. La lettre a été commencée le 7 septembre et continuée le 14. Elle est très longue (neuf pages dactylographiées). Il semble vraisemblable que l'envoi de Weil est arrivé entre le 7 et le 14.} (en français). La fureur, qu'avait anticipée Dieudonné, était bien présente:
\begin{quote}
Avez-vous l'idée d'un \og profiteur de guerre spirituel\fg? Il me semble que notre \og ami\fg André Weil soit un tel. Sûrement vous vous rappelez le problème principal de mes conférences à Paris en 1939, l'hypothèse de Riemann au champ des fonctions algébriques. Pendant que toutes les forces spirituelles de l'Europe étaient stérilisées par la guerre, Weil semble avoir conçu l'idée d'une démonstration de cette hypothèse-là. Sa première publication là-dessus, la note que vous m'aviez procurée, contient toutefois un lemme essentiel sans démonstration. Même Siegel et Chevalley n'en ont pu voir la validité. \`A leur demande directe de la démonstration du lemme, Weil a du [sic] concéder qu'il ne le puisse pas prouver. Or, il m'a fait parvenir une seconde note, dans laquelle il simplifie considérablement le contenu de la première, tout en se basant sur le même lemme, sans en donner la moindre indication d'une démonstration. Il dit seulement qu'une publication détaillée soit en \og préparation\fg. Quoiqu'il en soit, cette méthode me paraît tout autre que \og fair\fg. Supposons que Deuring ou moi, ayant beaucoup d'expérience dans ce champs-là, réussissions à prouver l'hypothèse de Riemann, soit en nous servant de la note de Weil, soit autrement. En ce cas, Weil, à coup sûr, prétendrait à la priorité, même si, à présent, il ne puisse point du tout démontrer son lemme. C'est ce que j'appelle une manière typiquement juive! Vous et moi, en des circonstances pareilles, n'auraient pas publié deux notes vanteuses.

Vous savez d'ailleurs que, quand j'étais sur le point de partir pour mes conférences à Paris en 1939, je reçus une lettre de Weil, dans laquelle il me donna des instructions détaillées à qui et dans quel sens j'eusse à m'adresser, pour lui obtenir une position à Paris\footnote{Il est question d'une lettre du 20 janvier 1939 (plusieurs mois avant le départ de Hasse pour Paris en mai) dans laquelle Weil demandait à Hasse, lorsqu'il serait à Paris, de lui rendre le grand service de suggérer aux professeurs parisiens de créer une chaire de théorie des nombres. Cette demande était un peu étrange, en effet, d'autant plus que Weil s'apprêtait à quitter la France pour échapper à la guerre qui s'annonçait (voir~\cite{AW-38}). L'\og insolence\fg supposée de Weil \og dès lors\fg serait plutôt, à la lecture de ses lettres, de la désinvolture. Il est vrai que sa dernière lettre se terminait par une injonction à Hasse de ne pas rédiger l'adresse de ses lettres en allemand gothique: \emph{In Strassburg können wenigstens die älteren Leute die deutsche Schrift noch lesen}, ce qui, au printemps 1939 et adressé à un Allemand, pouvait passer pour de l'insolence.}. Dès lors sa manière m'a paru comme le comble de l'insolence. Voici les lignes laconiques dont il accompagne le manuscrit de sa seconde note:

\emph{Dear Hasse, In the midst of the vastly more important affairs in which I hear you are at present engaged, you may still be able to spare a few minutes for the persual of the solution (\string!\string!\string!) of a problem you used to be interested in. With best greetings from the U.S.A. Yours sincerely .....}\footnote{Cette lettre de Weil ne semble pas se trouver dans le Nachla\ss\ Hasse. Peut-être ces quelques lignes étaient-elles écrites sur le tiré-à-part envoyé par Weil (\og manuscrit\fg ne semble pas être à prendre au pied de la lettre).}
\end{quote}

Rappelons que les \'Etats-Unis n'étaient pas encore en guerre et qu'il est possible que le courrier entre eux et l'Allemagne ait fonctionné à peu près normalement. Et remarquons que Weil a écrit à Hasse en anglais, alors que c'est en allemand qu'étaient rédigées toutes ses lettres des années 1930, mais Hasse était aussi bon anglophone que francophone\footnote{Sur la façon dont il perfectionna son anglais, avec Davenport, voir~\cite{Roquette04}.}.

\subsection*{Octobre 1941 (un article de Hasse)}
C'est le 20 octobre 1941 que le \emph{Jahresbericht} de la \textsc{dmv} reçut l'article~\cite{Hasse42} de Helmut Hasse. Il y faisait référence à l'exposé qu'il avait donné au cours de la réunion de la \textsc{dmv} de septembre 1938 à Baden-Baden: ainsi le \emph{Jahresbericht} était l'endroit naturel où publier ce texte. Il poursuivait:
\begin{quote}
Es war meine Absicht, diesen Vortrag, der zunächst nur einen skizzenhaften Überblick gab, bis in die Einzelheiten auszuarbeiten und als zusammenfassenden Bericht im Jahresbericht der Deutschen Mathematiker-Vereinigung zu veröffentlichen. Als ich im April 1940 zum Wehrmachtsdienst einberufen wurde, war die Arbeit an diesem Bericht, die immer weitere Kreise zog, kaum zur Hälfte gediehen, und seitdem bin ich nicht zur Fortführung gekommen. Ich habe mich jetzt entschlossen, das, was fertig vorliegt, zu veröffentlichen.\label{Wehrmachtsdienst}
\end{quote}

Il fait peu de doute que ce sont les deux notes de Weil et l'ironie qui avait accompagné l'envoi de la deuxième qui ont décidé Hasse à envoyer cet article. Ce n'est probablement pas par hasard qu'avril 1940 y est mentionné. L'article n'était pas prêt à être publié, les démonstrations n'étaient pas complètes. Certains des résultats furent démontrés par... van der Waerden (qui les publia plus tard dans~\cite{vdW47}). Hasse avait envoyé son article juste avant la réunion de la \textsc{dmv} qui se tint à Jena du 20 au 24 octobre 1941, au cours de laquelle il eut l'occasion de parler avec van der Waerden\footnote{Il n'est pas tout à fait hors sujet de signaler que, dans la recension qu'il fit pour \MR de l'article~\cite{vdW42} correspondant à l'exposé de van der Waerden à Jena, Chevalley lui reprocha d'avoir attribué à Kneser un résultat d'uniformisation locale que Kneser avait seulement annoncé, alors que Zariski l'avait démontré.} (l'un venait de Berlin, l'autre de Leipzig). Hasse lui envoya ensuite ses notes en lui demandant s'il pouvait démontrer certains des résultats qu'il y énonçait. \`A la fin de 1941, van der Waerden lui envoya ses démonstrations (sa lettre est datée de la Saint-Sylvestre), ce qui fut une notable exception à l'habituelle minceur de ses relations avec le groupe de Hasse.\label{demvdW} Pour davantage de détails, voir~\cite{schappa07}.

\subsection*{Novembre 1941 (recension de la note américaine de Weil par \MR)}
\MR publia une recension, toujours par Schilling, de la note américaine~\cite{Weil41}. Article publié en juillet, recension parue en novembre: sans doute un record de vitesse. Cette note ne fut jamais recensée, ni par le \Zbl, ni par le \JFM (qui a encore fonctionné pour les articles parus en 1941 et 1942).

\subsection*{L'année 1942  (recension de la note de Weil de 1940 dans le \JFM)}
Le fascicule du \emph{Jahresbericht} dans lequel est inclus l'article~\cite{Hasse42} de Hasse parut le 9 juin 1942. La date exacte à laquelle parut le \JFM de 1940 ne nous est pas connue avec précision, mais nous savons que ce fut en cette année 1942. 

\subsubsection*{Sur le \JFM}
Une précision s'impose ici. Nous sommes habitués à la recension rapide des journaux \Zbl et \MR. Mais ce fut une des  nouveautés de la création du \Zbl en 1931. Le \emph{Jahrbuch} de l'année $n$ est le volume qui recense tous les articles publiés pendant l'année $n$. D'ailleurs, dans ce volume, les références des articles ne contiennent pas l'année, qui est implicite. Le \JFM de l'année $n$ paraissait évidemment beaucoup plus tard; ces délais de publication (trois ans de retard en moyenne d'après~\cite{SiegmundS94}) devenus inadaptés à la production croissante furent même une des raisons de la création de \Zbl, qui publiait les recensions aussi vite que possible et était donc beaucoup plus rapide. Comme c'est indiqué sur sa couverture, le \JFM de 1940 date de 1942.

\medskip
Et, dans ce volume, parut donc la troisième recension de~\cite{weil40} par un algébriste allemand, beaucoup moins neutre que les deux précédentes. 
\begin{quote}
Verf. versucht, eine Lösung der Hauptprobleme der algebraischen Funktionenkörper mit endlich vielen Konstanten, insbesondere auch einen Beweis der Riemannschen Vermutung für Kongruenzzetafunktionen zu geben. Den Schlüssel dazu bieten einerseits die Severische algebraische Korrespondenztheorie, andererseits die Hurwitzsche transzendente Theorie der Korrespondenzen. Es handelt sich darum, die letzte Theorie in algebraischer Form auf die obigen Funktionenkörper zu übertragen. Dazu wird ein Isomorphismus der Divisorenklassengruppe nullten Grades und einer Gruppe einspaltiger Matrizen mit $2g$ Reihen modulo 1, deren Elemente einem Ring $\gotR$ angehören, hergestellt. Der Ring $\gotR$ entsteht aus dem Körper der rationalen Zahlen durch Komplettierung in bezug auf alle $l$-adischen Bewertungen, wo $l$ alle Primzahlen mit Ausnahme der Körpercharakteristik durchläuft. Nach dem Vorbild von Hurwitz wird jeder $(m_1, m_2)$-Korrespondenz eine quadratische, $2g$-reihige Matrix $L$ zugeordnet, deren Elemente kongruent 0 modulo 1 in $\gotR$ sind. Der inversen $(m_2, m_1)$-Korrespondenz entspricht ebenso eine Matrix $L^\prime$. Verf. formuliert nun einen entscheidenden Hilfs\-satz: Unter gewissen Bedingungen, die in der vorliegenden Note nicht angegeben sind, ist die Spur von $LL^\prime$ gleich $2m_2$. Aus diesem Hilfssatz folgt leicht die Riemannsche Vermutung. -- Die in dieser Skizze gegebenen Lösungsgedanken sind den deutschen Algebraikern nicht neu. Verf. hat die wesentliche Schwierigkeit in obigen Hilfssatz verlagert, dessen Beweis abzuwarten bleibt.
\end{quote}

Le troisième algébriste allemand qui recensait la note était Hermann Ludwig Schmid (1908--1956), qui avait passé sa thèse avec Hasse à Marburg en 1934, avait suivi à Hasse à Göttingen et y avait été son assistant deux ans (en 1935--37), était parti pour Giessen (où était Geppert), y était devenu en 1938 secrétaire éditorial de \Zbl, et se trouvait depuis 1940, à Berlin, l'assistant de Geppert (qui supervisait à la fois \Zbl et le \JFM)\footnote{Sur Schmid, voir les articles~\cite{JehneLamprechtBerlin,SiegmundSBerlin}.}. Signalons que le \JFM a republié, cette année-là, certaines recensions qui avaient été écrites pour et publiées par le \Zbl, et en particulier au moins deux de celles rédigées par van der Waerden. Dans le cas de la note de Weil, on a préféré une nouvelle recension. Il serait fort étonnant que les dernières phrases du texte de Schmid n'aient pas été inspirées par Hasse. Une hypothèse serait que, comme la publication de~\cite{Hasse42}, la parution de cette recension, alors que~\cite{weil40} avait déjà fait l'objet de celle de van der Waerden, était une réponse à l'ironie et aux fanfaronnades de Weil.

Notons encore que le caractère en principe \og plus objectif\fg~\cite{SiegmundS94} des recensions du \JFM est fortement compromis par la prise de position de celle-ci.

\medskip
Une autre remarque serait que Hasse ne se plaint nulle part, même pas dans les lettres à Julia, où pourtant il mentionne sa correspondance avec Weil, que les idées de sa lettre de juillet 1936 aient été utilisées (sans que celui-ci mentionne cette lettre) par Weil --- le peu d'étonnement que les algébristes allemands seraient en état de manifester vient pourtant aussi de là. La correspondance entre Hasse et Weil mériterait d'être publiée et commentée par un meilleur arithméticien que l'auteur de cet article. Le plus étonnant est que Weil, qui a, le 17 juillet 1936, accusé réception de la lettre de Hasse et discuté le fait que la théorie des correspondances pouvait être considérée du côté italien, semble avoir oublié ces questions jusqu'en février 1939. Lorsqu'il reprit la correspondance à se sujet, il avoua qu'il s'apercevait avec une grande honte (\emph{Mit grosser Beschämung}) que Hasse lui avait déjà envoyé une longue lettre en 1936...\label{weiletlalettre}

\subsection*{Octobre 1943 (recension de l'article de Hasse par \MR)}
Car ce n'est pas tout à fait terminé. L'acte suivant se joue dans \MR. Deux précisions. D'abord, les articles de Hasse ont été recensés, dans \MR pendant les années 40, par Alfred Brauer (mars 1940), Richard Brauer (février 1941, avril 1948, juillet-août 1948, avril 1949, juillet-août 1949, septembre 1949), recensions Brauer, donc... sauf deux, justement en 1943, où c'est Weil qui donna, pour le fascicule d'octobre, ses rapports sur les deux articles~\cite{Hasse42a,Hasse42}. Ensuite, s'il est vrai qu'André Weil a écrit beaucoup de recensions pour \MR, dès son arrivée aux \'Etats-Unis, il a surtout écrit sur la topologie\footnote{D'où son impressionnante connaissance de la bibliographie sur ce sujet, voir~\cite{CW-doc,AW-38}.}.

Lorsqu'il écrivit son rapport sur l'article de Hasse, Weil n'avait probablement pas encore vu la recension de Schmid, mais il fit entrer dans ses rapports la même ironie que dans la lettre à Hasse. Le premier article portait sur le théorème de Mordell-Weil. Voici la recension du deuxième, qui est, justement, celui du \emph{Jahresbericht}:
\begin{quote}
This, the author states, is the introductory section of a report, the writing of which was interrupted. The report was to cover the following topics, which the author lists as the most important results hitherto achieved in the study of algebraic curves over arithmetically significant fields of constants: (I) Weil's theorem of the finite basis [A. Weil, Acta Math. 52, 281--315 (1928)]; (II) Siegel's theorem on the finite number of integral solutions [C. L. Siegel, Abh. Preuss. Akad. Wiss. Phys.-Math. Kl. 1929, part 2]; (III) Lutz's theorem\footnote{Il s'agit d'\'Elizabeth Lutz (1914--2008), étudiante de Weil à Strasbourg et de l'article cité~\cite{Lutz37}, que Hasse avait publié dans le journal de Crelle. Cette publication était un des sujets de la correspondance entre Hasse et Weil en avril, mai 1936 et jusqu'à la lettre de juillet 1936.} on elliptic curves over p-adic fields [E. Lutz, J. Reine Angew. Math. 177, 238--247 (1937)]; (IV) Hasse's proof of the Riemann hypothesis for elliptic curves over Galois-fields [H. Hasse, J. Reine Angew. Math. 175, 55--62, 69--88, 193--208 (1936)] [this has now been superseded by the proof obtained for arbitrary genus by A. Weil, Proc. Nat. Acad. Sci. U. S. A. 27, 345--347 (1941); MR0004242 (2,345b)]; (V) Deuring's algebraic construction of the Abelian extensions of the function-field of an algebraic curve [M. Deuring, Math. Ann. 106, 77-102 (1932)]. These topics were to be treated from a ``unified point of view,'' which, needless to say, is supplied by the consideration of the Jacobian variety (or, in the language of the author, the Abelian function-field) belonging to the curve. The greater part of the present paper is devoted to an exposition of some of the more elementary properties of the Jacobian variety, in purely algebraic terms and without the use of the convenient tools supplied in the case of characteristic 0 by classical function-theory. These properties are of course couched in the arithmetico-algebraic language of the author and his school, which will be familiar to readers of his papers on elliptic function-fields, but may act as a deterrent on other classes of readers, and does not seem to the reviewer to be as well adapted to those questions as the language of algebraic geometry. Some new lemmas, mainly of technical interest, are stated and proved; some of the basic results of the classical theory are stated as unproved assertions in the abstract case. A final section is devoted to an exposition, restricted to algebraic curves, of Weil's so-called distribution-theory [loc. cit. chap. I], with the quantitative complements due to Siegel [loc. cit.]. One may express the hope that circumstances will soon restore the author to his mathematical studies and enable him to complete his report.
\end{quote}

Après la controverse mathématique (le langage arithmético-algébrique moins adapté que celui de la géométrie algébrique), la controverse personnelle. Si Hasse considérait que Weil avait \og profité\fg de la guerre, Weil ironisait sur le fait que cette même guerre empêchait Hasse de travailler... \og stérilisé par la guerre\fg, avait écrit ce dernier.

\subsection*{1946--48, 1951... fin de l'histoire?}
En 1946, la collection {Colloquium Publications} de l'\emph{American mathematical society}, publia le livre~\cite{Weil46} d'André Weil et, en 1948, c'est chez Hermann, dans la série \og Publications de l'Institut de mathématiques de l'université de Strasbourg\fg, que parut~\cite{Weil48a}. Ceci achevait la démonstration par André Weil de l'hypothèse de Riemann pour les courbes à corps de base fini

\begin{quote}
\noindent après huit années et plus de 500 pages, sa Note de 1940 est enfin justifiée!~\cite{SerreWeil}
\end{quote}

\enlargethispage{\baselineskip}

La guerre était finie, Weil avait de quoi triompher, l'hypothèse de Riemann et les résultats sur les fonctions $L$, même le lemme important de la note de 1940, tout était démontré, l'Allemagne était en ruines, Geppert s'était suicidé, le \JFM était vraiment mort, Schmid relançait le \Zbl, Weil ne parvenait pas à rentrer en France, Hasse traversait une sorte de (bref) purgatoire...

\medskip
En 1951, le \Zbl qui, sous la direction de Hermann Ludwig Schmid, avait recommencé à reparaître en 1947, publia une recension de~\cite{Weil48a}. Celle-ci contenait aussi implicitement une recension de~\cite{Weil46}, qui ne fut jamais traité comme tel par le \Zbl (pas plus que la note~\cite{Weil41} de 1941). Après les explications sur le caractère algébro-géométrique du travail et la positivité de la trace, les auteurs notaient en particulier (presque au détour d'une phrase):
\begin{quote}
\noindent Daraus wird die Richtigkeit des Analogons zur Riemannschen Vermutung in algebraischen Zahlenkörpern gefolgert.
\end{quote}

Ce résultat de Weil avait déjà été mentionné quelques mois plus tôt, par Hasse lui-même, dans sa recension de~\cite{Weil49} --- voir la note~\ref{noteweil49}. La recension dont il est question ici était elle aussi signée de Hasse (désormais professeur à Hambourg) et d'un de ses jeunes élèves, Peter Roquette (qui, né en 1927, était alors âgé de vingt-quatre ans).

\subsection*{\'Epilogue, épilogues?}
Outre les rivalités personnelles et les positions prises par rapport à la guerre, le débat portait aussi sur la nature des démonstrations, géométrie algébrique ou arithmétique? 

La démonstration de Weil, fondée sur la théorie des intersections et son livre \emph{Foundations of algebraic geometry}~\cite{Weil46}, était clairement dans le camp de la géométrie algébrique. Hasse et une partie des \og algébristes allemands\fg pensaient que la nature arithmétique de l'énoncé demandait une démonstration arithmétique. Comme nous l'avons signalé, le premier recenseur, van der Waerden, était d'ailleurs davantage du côté de la géométrie algébrique que ses collègues de l'entourage de Hasse (renvoyons à nouveau à~\cite{schappa07}). Et, en effet, une telle démonstration arithmétique apparut, dans le camp de Hasse, et c'est justement Roquette qui la donna, dans un article~\cite{Roquette53} dont le titre est \og Une preuve arithmétique...\fg. Sans surprise, ce furent Schilling dans \MR et Hasse dans le \Zbl (dès 1954) qui en écrivirent des recensions.

La géométrie algébrique avait pourtant fait ses preuves en arithmétique, rendant l'énoncé des conjectures de Weil~\cite{Weil49} possible\footnote{Puisqu'il est question ici de recensions... L'idée que l'on puisse calculer le nombre de solutions d'une équation (modulo $p$) par la formule de points fixes de Lefschetz, dont Serre~\cite{SerreWeil} dit \og Cette idée, vraiment révolutionnaire, a enthousiasmé les mathématiciens de l'époque\fg, ne semble pas avoir inspiré l'auteur de la recension de l'article~\cite{Weil49} par \MR, qui est passé complètement à côté du caractère révolutionnaire de l'article (et qui, cette fois, n'était pas Schilling). 

Mais cet article de douze pages a fait l'objet, dès mars 1950, d'une recension de deux pages et demie dans le \Zbl, par Helmut Hasse.\label{noteweil49}}: après le passage aux courbes de genre supérieur, la question était de passer aux variétés de dimension supérieure.

\medskip
Puisque nous en sommes aux \og héritiers\fg... signalons aussi, à la suite de~\cite{schappa06}, une \og petite phrase\fg, une note de bas de page dans la note historique aux chapitres I à~VII de l'\emph{Algèbre commutative} de Bourbaki (parue en 1965):
\begin{quote}
On sait que, malgré les efforts de Dedekind, Weber et Kronecker, le relâchement dans la conception de ce qui constituait une démonstration correcte, déjà sensible dans l'école allemande de Géométrie algébrique des années 1870--1880, ne devait que s'aggraver de plus en plus dans les travaux des géomètres français et surtout italiens des deux générations suivantes, qui, à la suite des géomètres allemands, et en développant leurs méthodes, s'attaquent à la théorie des surfaces algébriques: \og scandale\fg mainte fois dénoncé (surtout à partir de 1920) par les algébristes, mais que n'étaient pas sans justifier en une certaine mesure les brillants succès obtenus par les méthodes \og non rigoureuses\fg, contrastant avec le fait que, jusque vers 1940, les successeurs orthodoxes de Dedekind s'étaient révélés incapables de formuler avec assez de souplesse et de puissance les notions algébriques qui eussent permis de donner de ces résultats des démonstrations correctes.
\end{quote}

Les parenthèses et appositions, le changement de temps \og s'attaquent\fg au milieu d'une série de verbes au passé, cette trop longue phrase bien filandreuse peut difficilement être attribuée à Weil. \og Les notes historiques n'étaient jamais discutées dans les congrès Bourbaki\fg, m'écrit J-P. Serre. \og Les notes étaient surtout rédigées par Dieudonné --- mais avec la contribution de Weil quand le sujet en valait la peine. C'est sans doute le cas pour celle-ci (attention: parfois Samuel s'y intéressait aussi, et c'est peut-être le cas ici --- le caractère alambiqué de la phrase évoque pour moi Samuel plutôt que Weil).\fg Il est pourtant vrai que la fin, avec son imparfait du subjonctif et sa mention autrement incompréhensible de 1940, évoquent irrésistiblement la main de Weil.

\medskip
En 1979, dans le commentaire qu'il fit de la note~\cite{weil40} dans le premier volume de ses \OE uvres, André Weil, dont la preuve était complète depuis longtemps, et qui pouvait donc se le permettre, ironisa sur la recension de Schmid. Reprenons le paragraphe, dont des phrases ont été citées séparément plus haut, dans son ensemble:
\begin{quote}
En d'autres circonstances, une publication m'aurait paru bien prématurée. Mais, en avril 1940, pouvait-on se croire assuré du lendemain? Il me sembla que mes idées contenaient assez de substance pour ne pas mériter d'être en danger de se perdre. Si j'en crois Dieudonné, qui vit Hasse à Paris quelques mois plus tard, celui-ci s'en montra fort scandalisé. Le compte-rendu qui parut au \emph{Jahrbuch} de 1940 observa gravement qu'au \og lemme important\fg près, dont il \og convenait d'attendre la démonstration\fg, la note n'avait rien à apprendre \og aux algébristes allemands\fg. Faut-il en conclure que l'esprit de ceux-ci avait été quelque peu grisé par le succès de leurs généraux?~\cite[p.~546]{weilOC1}
\end{quote}

Si le commentaire sur le succès des généraux peut s'entendre comme un écho au dérapage antisémite de la lettre de Hasse à Julia (que certainement Weil n'a jamais vue), il semble un peu décalé par rapport à la réalité d'un Hasse empêtré dans les obligations militaires qui le stérilisent et l'empêchent de travailler comme il le souhaiterait.

\medskip
Hasse est mort en 1979 et n'a sans doute pas vu ce texte. C'est Peter Roquette qui écrivit une recension~\cite{RoquetteWeil} des \OE uvres de Weil pour le \emph{Jahresbericht} de la \textsc{dmv}.

\subsection*{Remerciements à...}

\begin{itemize}
\item Norbert Schappacher pour les documents qu'il m'a communiqués
\item Jean-Pierre Serre pour ses réponses à mes questions.
\item tous les lecteurs qui m'ont signalé des erreurs ou des fautes dans une version précédente, en particulier à Antoine Chambert-Loir, Nicholas Katz et Mathias Künzer.
\end{itemize}
\nocite{Beckeretal}

\newcommand{\SortNoop}[1]{}
\providecommand{\bysame}{\leavevmode ---\ }
\providecommand{\og}{``}
\providecommand{\fg}{''}
\providecommand{\smfandname}{\&}
\providecommand{\smfedsname}{\'eds.}
\providecommand{\smfedname}{\'ed.}
\providecommand{\smfmastersthesisname}{M\'emoire}
\providecommand{\smfphdthesisname}{Th\`ese}

\vfill

\begin{thebibliography}{10}

\bibitem{CW-doc}
{\scshape M.~Audin} -- \emph{{Correspondance entre Henri Cartan et André
  Weil}}, Documents mathématiques, Société mathématique de France, Paris, 2011.

\bibitem{AW-38}
\bysame , {\og {L'impossible retour Le cas d'André Weil}\fg},  (2011), en
  préparation.

\bibitem{Beckeretal}
{\scshape H.~Becker, H.-J. Dahms {\normalfont \smfandname} C.~Wegeler}
  (\smfedsname) -- \emph{{Die Universität Göttingen unter dem
  Nationalsozialismus}}, München, {K.G. Saur}, 1998, deuxième édition.

\bibitem{Deuring36}
{\scshape M.~Deuring} -- {\og {Automorphismen und Divisorenklassen der Ordnung
  $l$ in algebraischen Funktionenk\"orpern.}\fg}, \emph{Math. Ann.}
  \textbf{113} (1936), p.~208--215 (German).

\bibitem{Deuring37}
\bysame , {\og {Arithmetische Theorie der Korrespondenzen algebraischer
  Funktionenk\"orper. I.}\fg}, \emph{J. Reine Angew. Math.} \textbf{177}
  (1937), p.~161--191 (German).

\bibitem{Hasseoslo}
{\scshape H.~Hasse} -- {\og {\"Uber die Riemannsche Vermutung in
  Funktionenk\"orpern}\fg}, in \emph{{Congrès international des mathématiciens,
  volume 1}} (Oslo), A. W. Broggers Boktrykkeri, 1936, p.~189--206.

\bibitem{Hasse42a}
\bysame , {\og Der {$n$}-{T}eilungsk\"orper eines abstrakten elliptischen
  {F}unktionenk\"orpers als {K}lassenk\"orper, nebst {A}nwendung auf den
  {M}ordell-{W}eilschen {E}ndlichkeitssatz\fg}, \emph{Math. Z.} \textbf{48}
  (1942), p.~48--66.

\bibitem{Hasse42}
\bysame , {\og Zur arithmetischen {T}heorie der algebraischen
  {F}unktionenk\"orper\fg}, \emph{Jber. Deutsch. Math. Verein} \textbf{52}
  (1942), p.~1--48.

\bibitem{JehneLamprechtBerlin}
{\scshape W.~Jehne {\normalfont \smfandname} E.~Lamprecht} -- {\og Helmut
  {H}asse, {H}ermann {L}udwig {S}chmid and their students in {B}erlin\fg}, in
  \emph{Mathematics in {B}erlin}, Birkh\"auser, Berlin, 1998, p.~143--149.

\bibitem{Lutz37}
{\scshape E.~Lutz} -- {\og {Sur l'\'equation $y^2=x^3-Ax-B$ dans les corps
  $p$-adiques.}\fg}, \emph{J. Reine Angew. Math.} \textbf{177} (1937),
  p.~238--247.

\bibitem{Pekonen}
{\scshape O.~Pekonen} -- {\og L'affaire {W}eil \`a\ {H}elsinki en 1939\fg},
  \emph{Gaz. Math.} (1992), no.~52, p.~13--20, avec un commentaire d'Andr\'e
  Weil.

\bibitem{Roquette53}
{\scshape P.~Roquette} -- {\og Arithmetischer {B}eweis der {R}iemannschen
  {V}ermutung in {K}ongruenzfunktionenk\"orpern beliebigen {G}eschlechts\fg},
  \emph{J. Reine Angew. Math.} \textbf{191} (1953), p.~199--252.

\bibitem{RoquetteWeil}
\bysame , {\og {Weil, A. Collected papers}\fg}, \emph{Jahresbericht der
  \textsc{dmv}} \textbf{86} (1984), p.~29--32.

\bibitem{Roquette02}
\bysame , {\og The {R}iemann hypothesis in characteristic {$p$}, its origin and
  development. {I}. {T}he formation of the zeta-functions of {A}rtin and of
  {F}. {K}. {S}chmidt\fg}, \emph{Mitt. Math. Ges. Hamburg} \textbf{21} (2002),
  no.~2, p.~79--157, Hamburger Beitr{\"a}ge zur Geschichte der Mathematik.

\bibitem{Roquette04}
\bysame , {\og The {R}iemann hypothesis in characteristic {$p$}, its origin and
  development. {II}. {T}he first steps by {D}avenport and {H}asse\fg},
  \emph{Mitt. Math. Ges. Hamburg} \textbf{23} (2004), no.~2, p.~5--74.

\bibitem{Roquette06}
\bysame , {\og The {R}iemann hypothesis in characteristic {$p$}, its origin and
  development. {III}. {T}he elliptic case\fg}, \emph{Mitt. Math. Ges. Hamburg}
  \textbf{25} (2006), p.~103--176.

\bibitem{SchappaGottingen}
{\scshape N.~Schappacher} -- {\og {Der Mathematische Institut des Universität
  Göttingen}\fg}, in \emph{{Die Universität Göttingen unter dem
  Nationalsozialismus}}, 1998, {\rm in}~\cite{Beckeretal}, disponible
  sur~\url{http://www-irma.u-strasbg.fr/~schappa/NSch/Publications_files/MathI%
nstGoeNS.pdf}, p.~523--551.

\bibitem{schappa07}
\bysame , {\og A historical sketch of {B}. {L}. van der {W}aerden's work in
  algebraic geometry: 1926--1946\fg}, in \emph{Episodes in the history of
  modern algebra (1800--1950)}, Hist. Math., vol.~32, Amer. Math. Soc.,
  Providence, RI, 2007, p.~245--283.

\bibitem{schappa06}
\bysame , {\og {Seventy years ago: The Bourbaki congress at El Escorial and
  other mathematical (non)events of 1936}\fg}, \emph{Gac. R. Soc. Mat. Esp.}
  \textbf{11} (2008), no.~4, p.~721--735, Reprinted from Madr. Intelligencer
  {{\bf{2006}}}, 8--15.

\bibitem{SerreWeil}
{\scshape J.-P. Serre} -- {\og {La vie et l'\oe uvre d'André Weil}\fg},
  \emph{Enseign. Math.} \textbf{45} (1999), p.~5--16.

\bibitem{SiegmundSBerichter}
{\scshape R.~Siegmund-Schultze} -- \emph{{Mathematische Berichterstattung in
  Hitlerdeutschland}}, Studien zur Wissenschafts-, Sozial- und
  Bildungsgeschichte des Mathematik, vol.~9, Vanfenhoeck \&\ Ruprecht,
  G\"ottingen, 1993.

\bibitem{SiegmundS94}
\bysame , {\og ``{S}cientific control'' in mathematical reviewing and
  {G}erman--{U}.{S}.-{A}merican relations between the two world wars\fg},
  \emph{Historia Math.} \textbf{21} (1994), p.~306--329.

\bibitem{SiegmundSBerlin}
\bysame , {\og The {U}niversity of {B}erlin from reopening until 1953\fg}, in
  \emph{Mathematics in {B}erlin}, Birkh\"auser, Berlin, 1998, p.~137--141.

\bibitem{RSSvanderWaerden}
\bysame , {\og {Bartel Leendert van der Waerden (1903-1996) im Dritten Reich:
  Moderne Algebra im Dienste des Anti-Modernismus?}\fg}, in \emph{{Fremde
  Wissenschaftler im Dritten Reich Die Debye-Affäre im Kontext}}, Wallstein
  Verlag, Göttingen, à paraître.

\bibitem{vdW42}
{\scshape B.~van~der Waerden} -- {\og Die {B}edeutung des {B}ewertungsbegriffs
  f\"ur die algebraische {G}eometrie\fg}, \emph{Jber. Deutsch. Math. Verein.}
  \textbf{52} (1942), p.~161--172.

\bibitem{vdW47}
\bysame , {\og Divisorenklassen in algebraischen {F}unktionenk\"orpern\fg},
  \emph{Comment. Math. Helv.} \textbf{20} (1947), p.~68--109.

\bibitem{weil40}
{\scshape A.~Weil} -- {\og {Sur les fonctions alg\'ebriques \`a corps de
  constantes fini}\fg}, \emph{C. R. Acad. Sci. Paris} \textbf{210} (1940),
  p.~592--594.

\bibitem{Weil41}
\bysame , {\og On the {R}iemann hypothesis in function fields\fg}, \emph{Proc.
  Nat. Acad. Sci. U. S. A.} \textbf{27} (1941), p.~345--347.

\bibitem{Weil46}
\bysame , \emph{Foundations of {A}lgebraic {G}eometry}, American Mathematical
  Society Colloquium Publications, vol. 29, American Mathematical Society, New
  York, 1946.

\bibitem{Weil48a}
\bysame , \emph{Sur les courbes alg\'ebriques et les vari\'et\'es qui s'en
  d\'eduisent}, Actualit\'es Sci. Ind., no. 1041 = Publ. Inst. Math. Univ.
  Strasbourg {\bf 7} (1945), Hermann et Cie., Paris, 1948.

\bibitem{Weil49}
\bysame , {\og Numbers of solutions of equations in finite fields\fg},
  \emph{Bull. Amer. Math. Soc.} \textbf{55} (1949), p.~497--508.

\bibitem{weilOC1}
\bysame , \emph{{\OE uvres scientifiques, Volume~I}}, Springer, 1979.

\bibitem{weilsouvenirs}
\bysame , \emph{{Souvenirs d'apprentissage}}, Vita Mathematica, vol.~6,
  Birkh\"auser, Basel, 1991.

\end{thebibliography}
\end{document}